\documentclass{amsart}

\title{AP Theory IV: Intrinsic Topological Quantum Langlands Theory}

\author{H. E. Winkelnkemper}

\address{\begin{flushleft}\quad Department of Mathematics\\
\quad University of Maryland \\
\quad College Park, Maryland 20742
\end{flushleft}}

\email{hew@math.umd.edu}

\begin{document}
\maketitle

\begin{abstract}  Without using any moduli, sheaves, stacks, nor any analytic, nor category-type arguments, we exhibit an analogue to Geometric Langlands Theory in an entirely {\it model-independent}, non-perturbative,purely smooth topological context in Artin Presentation Theory. A basic initial feature is that AP Theory, as a whole, is already, {\it ab initio}, a universal canonical $2D$ $\sigma$-model, targeting smooth, compact, simply-connected $4$-manifolds with a connected boundary, and its topological Planckian quantum starting point, as well as its cone-like, {\it $\infty$-generated at each stage}, graded group of homology-preserving, but topology-changing transitions/interactions, exhibit the most general qualitative $S$-duality. We first point out the numerous mathematically rigorous, model-free, (i.e., intrinsic), topological AP analogues with the heuristic Kapustin-Witten version of Geometric Langlands theory, as well as the crucial differences between the two theories. The latter have to exist since AP Theory deals, a priori, essentially only with discrete group-theoretic {\it  presentation} theory, not Lie group representation theory, not with category nor  moduli theory, does not need classical SUSY, nor Feynman integrals, lattice models, etc., and furthermore is model-independent  and characteristic of real dimension four. It will become clear that AP Langlands theory is a model-independent, cone-like, graded, compact completion and topological 'envelope' of $4D$ $\mathcal N$=4 SUSY YM theory, their starting point, anchored, kept in place, so to speak, by a purely discrete group-theoretic analogue of Donaldson/Seiberg-Witten theory. Due to its discrete group-theoretic conceptual simplicity and universality, the AP Langlands program should also be considered to be a new type of intrinsic model-free symmetry, a universal Erlanger Program for Modern Physics. Thus in AP Theory, the analogue of the Geometric Langlands Program is essentially a model-independent super Erlanger Program, a topological Langlands program, maximally gauge-theoretic, with the maximum of discrete symmetries, i.e. crystallity, and the minimum of smooth topology, i.e., $2D$ flat membranicity.

\end{abstract}

\section{Introduction}

Due to the fact that AP Theory is a rigorous, but very conceptually simple theory, in order to understand the following, it is enough, at first, to just read the introductions and here indicated pages and sections of the papers by Witten, \cite{Wi}, Kapustin-Witten, \cite{KW}, Frenkel, \cite{F1}, \cite{F4}, and not all their technical analytic formulas, heuristic path integral and SUSY using or category-theoretic parts.

 For the very conceptually simple mathematically rigorous basics of AP Theory, see the short papers, \cite{W1}, p.6, \cite{W}. (We assume the reader has at least these six papers simultaneously at immediate hand, since at a first reading, the main point is just recognizing crucial undeniable analogies).

It is also useful, as a short overview of definitions and motivations, to read \cite{W3} and the introduction of  \cite{W4}, the first two sections of \cite{CW} and \cite{C}, as well as just the preface of E. Frenkel's book, \cite{F5}.

We interpret and exhibit the genuinely pure braid-theoretic AP Theory  as an intrinsic (i.e., model-independent) topological quantum Langslands theory (AP Langlands Theory, from now on) where all pertinent groups act as generally and freely as possible, \cite{W},section 2, and which is intrinsic, i.e., model-independent gauge theory-wise, just as Cobordism Theory is model-independent homology theory-wise, \cite{W1}, p.2.

{\it Thus already a basic point of this paper is: the gauge-theoretic Geometric Langlands philosophy of Kapustin-Witten, \cite{KW}, survives qualitatively, topologically, in this new very general model-independent, non-perturbative, highly symmetric, mathematically rigorous setting.} 

And, viceversa, the rigorous model-free AP Theory acquieres universal Langlands philosophy, so to speak, augmenting the rigorous qualitative physical arguments of \cite{W1}, \cite{W}.

Its existence, with all the properties of the title, is due mainly to the fact that, as explained in \cite{W}, p.5, AP Theory is the most conceptually simple, parameter-free, non-perturbative, background-independent, intrinsic, (i.e., model-independent), cone-like, graded  $(3+1)$-QFT, with a qualitative, group-theoretic Planckian $2D$ 'continuous to discrete' duality starting point,\cite{W}, p.2, a purely discrete group-theoretic fundamental equation, the Artin Equation, \cite{W1}, p.7, and with the most general, {\it $\infty$-generated at each stage}, cone-like, graded group of $Z$-homology-preserving, but topology changing transitions/interactions, our Torelli action,\cite{W}, p.5,  \cite{W1}, p.8, \cite{W3}.

In particular, smooth topology change in AP Theory, i.e., its Morse Theory, is purely group-theoretic and does not need 'sewing' nor 'bordism glueing' laws, etc., which are unnatural from the algebraic geometry point of view.

This sporadic $(3+1)$-QFT, characteristic of real dimension four, also has a complete Knot and Linking theory for all closed, orientable $3$-manifolds, see remark 5 in section 5 below.

We also point out, as stressed on pp. 3,5 of \cite{W}, that although in this sporadic $(3+1)$-QFT, the last surviving 'vestiges' after model-independence, of classical Hilbert space, are the symmetric integer matrices $A(r)$, \cite{W1}, p.9, the $\infty$ of the rigid static dimension of classical Hilbert space, a fundamental characteristic feature of classical QM, has only just been supplanted, in AP Theory, by the more dynamical (i.e., more physical) $\infty$ of the cone-like, {\it $\infty$-generated at each stage}, topology-changing Torelli action. One could say: classical Hilbert space has been 'quantized' {\it without} losing $\infty$ and that SUSY has  not {\it really} been 'broken', ($\infty$ is still there!), but has only been 'broken softly', (compare to \cite{BKLS}, p.9), graded in a mathematically rigorous, canonical, cone-like manner as stressed on p.6 of \cite{W}.

 {\it In AP Theory, in a non-perturbative manner, SUSY has just been membranically {\it tempered}, so to speak,to have the maximum of dynamics, namely cone-like,  $\infty$-generated at each stage, graded dynamics, act at most in dimension $4$ and is now completely mathematically rigorous.} Compare to  Witten, \cite{Wi7}, p.514, p.552. 

There is no need to introduce $\infty$-dimensional manifolds for this, compare to p.176 of Atiyah, \cite{A}.

Simultaneously, the construction of the smooth, connected, compact, simply-connected $4D$ manifolds $W^4(r)$, with connected boundary, from the discrete, i.e., $0D$ Artin presentation $r$, \cite{W1}, p.4, is, {\it in reverse}, the sharpest 'compactification' in the sense of \cite{KW}, \cite{BKLS}. 

{\it In AP theory we have the ultimate compactification, namely from $4D$ to $0D$, i.e.,  to discreteness.} 

This again can be considered to be a  'continous to discrete' Planckian phenomenon, a $4D$ one, caused by our fundamental topological $2D$ Planckian starting point, \cite{W}, p.2.

AP Langlands theory, a mathematically rigorous theory, which, as we explain in section 2 below, has a rigorous {\it qualitative} topological S-duality, is our extension of, and substitute for the heuristic maximally supersymmetric extension of YM theory, i.e., $4D$ $\mathcal N$=4 SUSY YM theory, with its classical SUSY and only {\it conjectured} Montonen-Olive EM duality, the starting point of Kapustin and Witten's version of Geometric Langlands theory.

 This extension is model-independent and does not depend on any particular Lie group $G$, i.e., is intrinsic in this manner also, and is kept in place, so to speak, by containing a purely discrete group-theoretic version of Donaldson/Seiberg-Witten theory, despite the fact that, a priori, moduli do not exist in AP Theory, see \cite{W1}, p.9, \cite{CW}, \cite{R}, p.621, and compare to \cite{KW}, p.6 , \cite{Wi3}, p.4.

The best way to first  understand the purpose of this paper is to understand and realize that it is an intrinsic, meta-mathematically independent incarnation of the 'topological sector' of Geometric Langlands as mentioned on p. 13 of Frenkel's Bourbaki article, \cite{F1}. Compare also to \cite{FLN}. 

In a nutshell, we do for the model-free gauge theory approach to Geometric Langlands \'a la \cite{KW}, (see also \cite{Wi4}, p.1 and section 3.5), what we did for QFT in our previous paper, \cite{W}, p.5.  We  expand to the purely discrete group-theoretic maximum, the 'topological sector' on p.13 of \cite{F1}, with the sporadic, characteristic of dimension $4$, model-independent, cone-like, graded $(3+1)$-QFT of \cite{W}, p.5, instead of other more coarse, 'bordism glueing' using TQFTs, or classic heuristic SUSY using theories. 

Then without any need for any elaborate gauge theoretic machinery in order to "twist" $4D$ $\mathcal N$=4 Super YM theory, etc., AP Theory reveals itself, {\it ab initio}, as a natural, universal, canonical $2D$ $\sigma$-model targeting smooth, compact, simply-connected $4D$ manifolds with a connected boundary. See section 2 below.

The more dynamic, cone-like, graded,  $\infty$-generated at each stage, SUSY of \cite{W}, p.6,  substitutes for classical SUSY. It is also very relevant that the basic equation of AP Theory , the Artin Equation, is a purely discrete group-theoretic equation, \cite{W1}, p.7, \cite{W}, p.2, thus avoiding all types of analytic and nomographical obstructions. For example, there is no need in AP Theory to 'place some supersymmetric boundary condition' as on p.112 of \cite{KW}, nor have to worry about it, see \cite{GaWi}.

The ingredients of AP Langlands theory are simply $r$, $h(r):\Omega_n\to \Omega_n$, (where $\Omega_n$ denotes the compact, smooth, genus zero, flat $2D$ disk with $n$ holes) and the $W^4(r)$, see \cite{W}, p.2, \cite{CW}, p.65, all acted upon naturally, in unison, by the cone-like, graded, $\infty$-generated at each stage, homology-preserving, but topology-changing Torelli transitions/interactions, which now serve as $S$-operators, see \cite{W}, p.3 and compare to \cite{KW}, p.74.

In this paper these transitions/interactions should be considered as the most general set of $S$-operators, see the next section.

When the Torelli act on the $h(r):\Omega_n\to \Omega_n$ they are the analogue of the 'rigid surface operators' of \cite{Wi}, \cite{KW}, \cite{GW}, \cite{KS}, \cite{KSV}.

 In AP Theory, roughly put, the  Geometric Langlands correspondence reduces to the construction of the $4D$ branes $W^4(r)$, from the discrete Artin presentation $r$, via the $2D$ branes, $h(r):\Omega_n\to \Omega_n$, all acted upon,  {\it in unison}, by the cone-like, $\infty$-generated at each stage, graded analogue of SUSY described above, compare to \cite{Wi}, p.33, \cite{K2}. It is the most dynamic, universal, intrinsic, topological quantum theory 'pushing' the Langlands correspondence and should be considered as the meta-mathematical limit, the most general, model-free manifestation of Langlands philosophy, augmenting and justifying its Erlanger Program similarity, see \cite{Ge}, p.178.

Thus AP Langlands essentially justifies Witten's 'third motivation' of p.6 of \cite{Wi}.

 The smooth $4$-manifolds $W^4(r)$, our 4-dimensional branes, are actually constructed from the  $2$-dimensional Ur-branes $h(r):\Omega_n\to \Omega_n$, see \cite{W1}, section 2, and thus AP Theory, as far as topology is concerned, is essentially a purely $2D$ membranic theory: the smooth topology of these 'membratons' is all that is needed, as far as even the $4D$ smooth topology is concerned, due to the fundamental construction, a $4D$-brane construction, of AP Theory, see \cite{W1},p.4. Thus AP Theory is the most 'quantic' of all membrane-using theories, but which is dynamically globalized by the above cone-like, $\infty$-generated at each stage, graded group of Torelli transitions/interactions.

Our theory is just  the 'cobordic', model-free, non-analytic version of Kapustin-Witten Theory, with their electro-magnetic duality in their $4D$ SUSY gauge theory substituted by the membranic, 'topological EM theory' of \cite{W}. This topological theory has physical connotations: it can be considered to be the membranic equivalent of the Strong Force, since the $h(r):\Omega_n\to \Omega_n$, as equivalence classes modulo boundary-fixing isotopies,should be considered membranic 'hadrons', since they confine their quarks, partons , as explained on pp.4, 11, of \cite{W}, \cite{W1}, p.11, thus giving a rigorous topological analogue of 'confinement'. Compare to, e.g., \cite{H}, \cite{HSZ}, as well as Wilczek, \cite{Wilc}, p.85.

It is thus a mostly discrete  Dynamic Theory, which exploits the $2D$ genus zero 'flatness' of the 'membratons', 'hadrons' $h(r):\Omega_n\to \Omega_n$, \cite{W}, p.3, via the purely discrete group-theoretic Artin presentation $r$, with the Torelli action, to the maximum, in a clear, universal, conceptually simple Erlanger Program-like discrete group-theoretic style. The discrete Artin presentation, $r$, so to speak,in a subtle manner, gives a "quantization of the flat {\it membranic} topological 'connection' $h(r)$" on the compact 2-disk $\Omega_n$, thereby solidifying, crystallizing, quantizing its 'membranicity' and planarity, compare to 't Hooft, \cite{H}, \cite{W},p.2. For all this it is necessary that the topological genus of $\Omega_n$ is zero.

It can also be viewed as the primordial relation of the discrete 'string' $r$, \cite{W}, p.7, with the membrane $h(r):\Omega_n\to \Omega_n$. Our discrete, (i.e., lowest-dimensional) 'strings', the Artin presentations, $r$, (\cite{W}, p.7), are related to $2D$ flat membranes, ab initio, from the beginning.

This membranic topological  flatness, planarity, is all that remains in AP Langlands of flat connections on Riemann surfaces, which have the {\it 'least energy, most SUSY'}, etc., see \cite{Wi}, p. 56.

We also point out that, just as at the end of \cite{W}, the huge class of smooth topology changing, but homology preserving Torelli transitions/interactions, now operating as $S$-operators, should fatally obstruct any rigorous derivation  of some of the analytic desiderata of \cite{KW}. Compare to \cite{W}, p.12, especially as their operator theory is concerned. For example, in AP Theory operators can change the smooth structure of a $(3+1)$-manifold, in a purely discrete group-theoretic manner, but leave the underlying topological structure intact, see \cite{W}, p.5, \cite{C}. See also remark 1 in section 6 below.

 The sheer existence of AP Langlands theory seems to answer a fundamental question of Geometric Langlands Theory, Ben Zvi's " What to do with the trivial local system?". See {\it 'Generalized Geometric Langlands is false'} in {\it The n-category Caf\'e} blog.

AP Langlands Theory is the universal envelope of Geometric Langlands from the {\it intrinsic} model independent, non-perturbative gauge theory point of view; compare to \cite{Wi4}, p.1, \cite{Wi2}.

The gauge-theoretic point of view has been maximized in a non-perturbative, parameter-free, background and model-independent manner. 

In fact it seems to incarnate a Super Erlanger Program for Modern Physics, with a very, very large amount of 'bonus symmetries', compare to \cite{I}, thus shedding new light, via its new infusion of pure discrete group theory, into $4D$ smooth manifold theory, on very important, but stagnant, open problems, such as that of the cosmological constant, dynamic dark energy problem, the Clay YM  Mass Gap problem and QCD Confinement problem, see \cite{W}. 

This and its model-independence (as in Cobordism Theory) are the main reason that AP Langlands theory, one could say, AP Theory interpreted as a Langlands Theory,is immediately related, by undeniable analogies to YM Theory, QFT, QCD, Holography, etc.,see \cite{W1}, \cite{W}.

AP Theory is a discrete purely group-theoretic theory except only for the adjoining of the smooth, $2D$ hadrons, the flat $2D$ membranic $h(r):\Omega_n\to \Omega_n$. All $4D$, $3D$, topology now follows, without any new analytic, nor topological axioms; smooth topology has been reduced to the minimum and hence pure discrete group theory to the maximum. 

In AP Theory,in an {\it intrinsic, model-independent} manner, discrete gauge-theoretic symmetry has been maximized, and topology has been minimized, exhibiting it as a genuine super dynamical gauge-theoretic Erlanger Program for Modern Physics.

\section{AP Langlands Theory and Kapustin-Witten's Geometric Langlands Theory}

The main pertinent, basic mathematically rigorous topological facts and analogues in AP Theory to their theory (compare to \cite{Wi}, pp.60-75, and their main ideas in the introduction of \cite{KW}), are:

I. AP Theory already, ab initio, is a natural, intrinsic, {\it canonical} $2D$ $\sigma$-model, where the target manifolds of the discrete, topology-less Artin presentations $r$, via the $2D$ flat membranic hadrons  $h(r):\Omega_n\to \Omega_n$, are the compact, smooth, simply-connected $4$-manifolds with connected boundary,, our $4$-branes, the $W^4(r)$, \cite{W1},p.4, \cite{W}, p.2,  instead of Hitchin's moduli spaces; compare to \cite{KW}, p.3, \cite{F1}, p.14.

Thus there is no need of the 'twisting procedure' of \cite{Wi}, p.60, \cite{Wi3}, p.4, \cite{KW}, p.58, on $4D$  $\mathcal N$=4 Super YM theory, to a $2D$ $\sigma$-model, the starting point of \cite{KW}. It already is, {\it ab initio},  a mathematical 'given', in AP Langlands theory.

The sporadic $4D$-brane construction of the $W^4(r)$ from the $2D$ $h(r):\Omega_n\to \Omega_n$ is already the universal all encompassing AP $\sigma$-model, where all the following takes place, where literally 'all the action' takes place.

{\it In AP Langlands, the $\sigma$-model analogy never breaks down.}

II. Our starting 'electric-magnetic duality' is the AP theory smooth topological Planckian analogue to Planck's postulate, see p.2 of \cite{W}, intuited already, it seems,  by Poincar\'e himself with his '{\it resonateurs}', see \cite{W}, p.12. That is our {\it rigorous} topological quantum S-duality, which is globalized dynamically by the above Torelli action acting in unison on the discrete, i.e.,$0D$ $r$, the $2D$ $h(r):\Omega_n\to \Omega_n$ and the $4D$ $W^4(r)$.

In AP Theory we have a basic crystallic- membranic  duality between the topology-less, i.e., zero-dimensional, crystallic,  'electric' Artin presentation $r$ and the isotopy class of diffeomorphisms, the $2D$ flat membranic 'magnetic' $h(r):\Omega_n\to \Omega_n$, so to speak. 

As mentioned already this also has meaning as the immediate first relation of the discrete {\it 'quantum string'} $r$, see \cite{W}, p.7, with $2D$ flat {\it membranes}, namely $h(r):\Omega_n\to \Omega_n$.

This is the rigorous 'local'{\it topological} analogue to the still only conjectured classical S-duality as mirror symmetry of \cite{KW}'s $\sigma$-models. Here it is relevant to point out that all Artin presentations on two generators are in one to one correspondence with all $2\times 2$ symmetric integer matrices, \cite{W3}, \cite{R}, p.619, the last {\it qualitative} vestige, after model-independence, of classical quantitative $S$-duality.

 There is a subtle and important point here: the quantitative {\it conjectural} S-duality \'a la Olive-Montonen, as on p.2 of \cite{F1}, is supposed to be exactly one-to-one, when using $4D$ $\mathcal N$=4 SUSY gauge theory. Our rigorous topological quantum S-duality, between $r$ and $h(r):\Omega_n\to \Omega_n$, though, as explained on pp. 4, 11, of \cite{W}, is  exactly one-to-one only if we consider the diffeomorphism  $h(r):\Omega_n\to \Omega_n$ , \cite{W}, p.3, \cite{W4}, p.226, determined by $r$, only up to smooth boundary-fixing isotopies of $\Omega_n$, \cite{W}, p.3, \cite{W4}, p.226. These equivalence classes, our membranic hadrons, now give  an analogue of 'confinement' in AP Theory, as explained on p.4 and p.11 of \cite{W}.
 Thus, ab initio, in AP Langlands theory, S-duality is related to the important concept of {\it confinement} of QCD, in a rigorous manner, \cite{W}, pp.4, 11.

This rigorous local $S$-duality leads naturally to our global $S$-operators, the Torelli.

Thus $S$-duality in AP Theory is not quantitative, but qualitative and dynamic, similar to our versions of the cosmological constant, dynamic dark energy and the  YM Millenium mass gap problems in \cite{W}.

III. The sheer existence of AP Theory should be considered the meta-mathematical 'limit', so to speak, when the closed Riemann surface $X$ of p.14 of \cite{F1}, or $C$ of \cite{Wi}, p.32,p.65, \cite{KW}, pp. 30, 72, 112, \cite{Wi4}, p.10, 'shrinks' to a point. 

This is the mathematically rigorous most radical meta-mathematical 'limit' of \cite{Wi}, p.32, p. 65, \cite{KW}'s {\it ' $\Sigma$ large and  $C$ small'}. See also \cite{F1}, p.14.

 This, our  {\it compactification}, in the sense of \cite{KW}, p.3, is as radical as possible, from $4D$ to $0D$, i.e., discreteness, but still mathematically rigorous, and represented by a canonical, universal $2D$ $\sigma$-model, targeting smooth, compact, simply-connected $4$-manifolds with a connected boundary: the whole of AP Theory.

 The targets of  the topological, smooth, membranic hadrons $h(r):\Omega_n\to \Omega_n$, (our most 'local' topological concept),see \cite{W}, p.3, are the $4D$ 'spacetime universes', the  $4D$-branes, $W^4(r)$, (our most 'global' concept), instead of the  Hitchin moduli manifolds of \cite{KW}, p.3. Their Riemann surface $X$ has been rigorously  'planarized', into a hadron $h(r):\Omega_n\to \Omega_n$ in the 'limit', {\it but} still leaving a non-trivial model-free, non-perturbative Langlandsian theory in its wake; compare also to \cite{Wi5}, p.1.

This hadron has 'minimum energy and maximum symmetry', more so than Witten's flat bundles over Riemann surfaces, \cite {Wi}, p. 56.

Due to its qualitative universality and the starting point with the $h(r):\Omega_n\to \Omega_n$, there is no distinction between 'un-ramified' and 'ramified' versions of AP Langlands theory, compare to \cite{GW1}, \cite{F3}, p.2.

We can say: AP Theory 'envelopes' $4D$ $\mathcal N$=4 SUSY YM theory {\it membranically} and so 'tightly' that even a very non-trivial, purely group-theoretic Donaldson/Seiberg-Witten theory is still in place, \cite{CW}, \cite{W1}, p.9.

This goes beyond the mere $2D$ $\sigma$-model picture and should have important applications elsewhere.

This radical but still  mathematically {\it rigorous}  compactification, makes AP Theory a very conceptually simple Langland-sian theory with  clear, rigorous, topological qualitative analogues to the main features of \cite{KW}, and others.

This completely justifies our $(3+1)$-QFT of section 2 of \cite{W}, as the correct rigorous 'enveloping' analogue to $4D$ $\mathcal N$=4 SUSY YM theory;  $r$ is electric, $h(r):\Omega_n\to \Omega_n$ is magnetic, crystallic is electric, membranic is magnetic.

This is the AP Theory non-perturbative, model-independent, topological strong force analogue to the EM of Montonen-Olive in $4D$ $\mathcal N$=4 SUSY YM theory, compare to \cite{Wilc}, p.85. 

 It also re-inforces and justifies the terms, Gauge Theory and  YM in \cite{W}, \cite{W1}, even in these model-free cases.

IV. All operators in our theory are simply part of the discrete group-theory of AP Theory and they are automatically 'rigid' in the sense of \cite{GW}, since AP Theory is a parameter-free theory. Their type is only distinguished by their group-theoretic properties and/or  by {\it where}, in which dimension, we make them operate. Compare to \cite{KW}, p.85. In AP Theory operators are all embedded in rigorous classical pure braid group theory, where they act in unison on the discrete,i.e., $0D$ $r$, the $2D$ $h(r):\Omega_n\to \Omega_n$ and the $4D$ $W^4(r)$, i.e., on our whole $\sigma$-model, exhibiting AP Langlands theory, in particular, as a universal $\sigma$-model with a super Erlanger Program; compare to \cite{Ge}, p.178.

Thus in AP Langlands theory the operators, as simply elements of the cone-like, graded Torelli group, are acting on the fundamental $\sigma$-model itself, (compare to the 'operation S' of \cite{KW}, p.74), naturally in unison, on the '0-branes', the discrete $r$ as well as on the 2-branes $h(r):\Omega_n\to \Omega_n$ and the 4-branes $W^4(r)$, compare to p.74 of \cite{KW}. In that order, they should be called {\it rigid} string operators, {\it rigid} surface operators and {\it rigid} Hecke operators, compare to the line operators in \cite{KW} and surface operators in \cite{GW}.

All operators in AP Langlands theory are analogues of 'S-operators'.

It is the Torelli group which represents the {\it Geometry}, as well as the Algebra, of the Wilson, 't Hooft and Hecke operators in a truly Erlanger Program group-theoretic style. Compare to \cite{KS}.

Due to their relationship with confinement, (see \cite{Wi1}, p.2, \cite{W}, p.4, section 5), it is natural to consider our operators as 'Wilson operators' when acting on the $h(r):\Omega_n\to \Omega_n$, as ' 't Hooft operators' when acting on the discrete topology-less Artin presentations $r$. 

Because of the following, the Torelli acting on the 4-branes $W^4(r)$, are the AP-analogues of Hecke operators, the $4D$-branes $W^4(r)$ representing intrinsic 'Hecke eigensheaves', compare to \cite{Wi}, p.63.

In AP Langlands theory,'Hecke modifications', \cite{KW}, p.135, are caused by the topology change of the Torelli, acting on the $4$-branes $W^4(r)$ and their connected boundaries.
In AP Langlands theory,topology change is the genuine analogue to Hecke modification, since the latter uses glueing methods  and there is no  glueing in AP Theory. This logically forces our Torelli, when acting on the $W^4(r)$, to be the analogues of Hecke and $S$-operators in AP Theory.
 
There is no need for difficult procedures of inserting singularities, using modifications, nor path integral methods in operators in  Quantum AP Langlands theory, compare to p.5 of \cite{Wi3}, \cite{KW}, p.120.

Here perhaps it is instructive to compare the $h(r):\Omega_n\to \Omega_n$ to the figure of p.16 of \cite{Wi4}. They represent 'topologically flat, planar connections', so to speak,  on the multiply  'punctured', but compact 2-disk, $\Omega_n$. They are our topological analogues to the 'opers' of p.vi of \cite{F5}.

Due to the universality of AP Theory, the Artin presentations $r$ can be called 'electric eigenbranes', compare to \cite{Wi}, p.62, their intrinsic 'eigen property', simply being the preservation, under the Torelli action, of the matrix $A(r)$ obtained from $r$ by abelianization, see \cite{W1}, p.7, which determines the $Z$-homology of $W^4(r)$ and its always connected boundary.
Similarly, in $2D$, the $h(r):\Omega_n\to \Omega_n$ are 'magnetic eigenbranes', with their intrinsic model-free 'eigen property' now  being that they induce the identity on the homology of $\Omega_n$, see \cite{W4}, p.244.

The radical intrinsic, universal, (i.e. model-independent) compactification of AP Langlands Theory, forces our Hecke 'eigen property' to be intrinsic, model-free, 'exterior', i.e., what is left invariant  under the Torelli operators is the $Z$-homology of the 4-manifolds and their boundaries as explained in \cite{W}, \cite{W1}, p.8. This also is a manifestation of model-free 'eigen' rigidity.

  Thus, although, a priori, in AP Langlands, due to its model-free intrinsicity, {\it all} branes are eigen-branes, their 'eigen property' is that they preserve the $Z$-homology of the $W^4(r)$ and their boundaries, see \cite{W1}, p.8, and their 'modification' theory is just their topology-change theory, see \cite{W1},p.8.

In the model-independent AP Theory, the absolute 'eigen' property, is that all pertinent {\it $Z$-homology} be preserved, under the radical topology-changing Torelli action. It is the last vestige of 'eigen-ness' after the model-independence caused by the Torelli action.

Questions of 'periodicity' as in \cite{KW} should be discussed later, as, e.g., periodicity under the Torelli action, or better adressed as questions about the periodic points of the $h(r):\Omega_n\to \Omega_n$, our 'quarks', 'partons', see \cite{W1}, p.11, which the membranic hadrons $h(r):\Omega_n\to \Omega_n$ hold together. 

{\it This leads naturally to zeta-functions, the Artin-Mazur zeta-functions of the diffeomorphisms $h(r):\Omega_n\to \Omega_n$, thus giving another important Langlands property}.
 
In AP Langlands all operators live group-theoretically in the same graded group, act in unison, preserving the corresponding $Z$-homology, and are distinguished only by names, or group-theoretically by whether they operate on the discrete $0D$ $r$, the $2D$ $h(r):\Omega_n\to \Omega_n$ or the $4D$ $W^4(r)$.

V. {\it Now}, the  AP-Langlands analogue to the classical Geometric Langlands Duality : "Flat $G^L$-bundles on $X$ and D-modules on $Bun_G$" of, e.g., p.71 of \cite{F4}, \cite{KW}, pp. 122, 184, 198, is simply our canonical $\sigma$-model, where the membranic hadrons, (our $2D$ Ur-branes), $h(r):\Omega_n\to \Omega_n$ target the smooth, compact, simply-connected $4D$ manifolds $W^4(r)$, our $4D$-branes. The latter are the analogue of '$\mathcal D$-modules on $Bun_G$', now, {\it as in classical string theory}, in function of  D-branes, as in \cite{Wi}, p.63,  \cite{F4}, p.71, \cite{KW}, p.112, 184, 198, and the former, on the left Galois side, see \cite{F2}, p.156, p.166, are the analogues of 'flat $G^L$-bundles on $X$', equivalently, of homomorphisms of $\pi_1(X)$ to $GL_n(C)$. 

Instead of flat bundles on the Riemann surface $X$, i.e.,representations of its fundamental group into $GL_n(C)$, we have, due to our compactification of III. above, the more canonical Artin presentation $r$, which uniquely determines the isotopy class of intuitively 'flat' diffeomorphism $h(r):\Omega_n\to \Omega_n$. \cite{W}, p.2. These diffeomorphisms then determine the $4D$ branes $W^4(r)$, as explained on p.4 of \cite{W1}. These are our absolute, intrinsic, topological analogues to the 'local geometric Langlands parameters' and 'loop groups' (see \cite{F5}), are substituted by the $\pi(r)$, the fundamental groups of all connected, closed, orientable 3-manifolds.

 As mentioned above, this is the AP answer to Ben-Zvi's fundamental question about geometric Langlands: "What to do with the trivial local system?". 

 In AP Theory, the $3$-manifolds, $M^3(r)$, (the connected boundaries of the $W^4(r)$), and certain {\it presentations} of their fundamental group, $\pi(r)$, (namely Artin presentations), replace closed complex Riemann surfaces $X$ and certain {\it representations} of $\pi_1(X)$, (namely representations into $GL_n(C)$).

As concerning the right 'automorphic' side, the $4D$ branes $W^4(r)$, with their topology-changing  Torelli transitions are all that is left after, instead of using the 'holomorphic sector of a CFT', see \cite{Wi6},p.1, (where this idea seems to have first started), we expand our 'topological sector' to the maximum with the sporadic $(3+1)$-QFT as explained in the introduction above.

In AP Theory, Langlands duality boils down to the Planckian $2D$ 'continous to discrete' starting point, \cite{W}, p.2, and the fundamental smooth $4D$ brane construction ot the $W^4(r)$, via the $h(r):\Omega_n\to \Omega_n$ from the discrete Artin presentation $r$, completed dynamically with the {\it $\infty$-generated at each stage}, graded Torelli action. 

Due to its intimate relation with pure framed (i.e. colored)  braid theory, \cite{W1}, p.7,  we can also invoke Galois Theory in a different manner, via Neukirch's paper relating Galois groups to braids, \cite{N}, as well as via \cite{Ih}.

Notice also that colored, (i.e. framed) braids, (which are in one-to-one correspondence with Artin presentations,\cite{W1}, p.7), are mentioned by Langlands himself in \cite{L}, p.13, p.18, p.19. 

Instead of a specific quantitative classical $S$-duality in $4D$ $\mathcal N$=4 Super YM theory, being the analogue of Langlands duality, as in the \cite{KW} theory, this duality is now represented in AP Theory,{\it qualitatively and dynamically} by the cone-like, $\infty$-generated at each stage, graded group of homology-preserving, but smooth topology-changing Torelli transitions/interactions, see \cite{W}, p.5, which, furthermore, acts naturally in unison on the $r$, the $h(r):\Omega_n\to \Omega_n$ and the $W^4(r)$. Compare to \cite{Wi}, p.33. They incarnate  at the same time 'mirror' symmetry for our canonical $\sigma$-model the whole of AP Theory, compare to p. 60 of \cite{Wi}.

All these constructions are natural with respect the powerful $\infty$-generated at each stage, graded topology-changing Torelli transitions/interactions. The latter transitions/interactions should be considered the analogues of the 'ordinary differential operators' of \cite{KW}, p.198, and, at the same time, the AP $S$-operators, see \cite{KW}, p.74.

Thus the left side of AP Langlands starts at the most basic topological Planckian quantum level in the sense of \cite{W}, p.2. Compare to \cite{Wi3}, \cite{Wi5}. This Planckian duality is our rigorous topological quantum S-duality, our crystallic to membranic duality, which when dynamically globalized by the Torelli action, is  our analogue to the conjectured Montonen-Olive electro-magnetic duality of $4D$ $\mathcal N$=4 SUSY YM theory, the starting point of \cite{KW}.

In AP Theory, the quantitative, analytic \cite{KW} Langlands philosophy boils down {\it qualitatively, topologically}, in a model-free manner, to: AP Theory is a canonical universal  $2D$ $\sigma$-model, a natural compactification to $2D$ of the $4D$ manifolds $W^4(r)$, where both sides are determined by the $0D$, i.e. discrete Artin presentation $r$, and the group of $S$-operators acting on it are the Torelli above. 

Thus in AP Theory, Geometric Langlands Duality essentially reduces {\it qualitatively} to the sheer existence of the cone-like, graded, membranically tempered SUSY of \cite{W}, p.6.

Among all the Langlandsian theories, AP Langlands theory is the compact, universal topological envelope theory of Geometric Langlands Theory from the point of discrete group theory, i.e.,intrinsic, non-perturbative gauge theory, compare to \cite{Wi4}, p.1. In fact, in AP Theory the analogue of the geometric Langlands Program actually  becomes a Super Erlanger Program for Physics, with the maximum of discrete group theory and the minimum of smooth topology, namely $2D$ flat smooth membrane topology.
 
AP Langlands Theory, by making the gauge theory model-free (just as Cobordism Theory made homology theory model-free) envelopes the model-dependent maximally supersymmetric extension of YM-theory, namely $4D$ $\mathcal N$=4 SUSY YM Theory, the starting point of \cite{KW}, and this 'enveloping' is so 'tight', that a purely discrete group-theoretic version of Donaldson/ Seiberg-Witten theory remains, see \cite{W1}, p.9, \cite{CW}, p.2, and compare to \cite{Wi2}.

 {\it Thus AP Langlands Theory and Kapustin-Witten's Geometric Langlands theory are rigorously related meta-mathematically in the most group-theoretic, i.e., gauge-theoretic fashion.}

 We can say AP Langlands theory is the model-independent intrinsic topological 'completion', the  meta-mathematical 'halo', so to speak, of Kapustin-Witten's Geometric Langlands Theory, anchored and kept in place by a purely group theoretic Donaldson/Seiberg-Witten theory, \cite{W1}, p.9.

 Despite its model-independence, AP Theory envelopes the maximally supersymmetric extension of Yang-Mills theory, namely $4D$ $\mathcal N$=4 Super YM theory in a very meta-mathematically tight manner, with Donaldson/Seiberg-Witten theory, as well as Kapustin-Witten theory, still well represented.

In fact, the purely group-theoretic  Donaldson Theorem of \cite{W4}, p.240, \cite{W1}, p.9, \cite{R}, p.621, relating number-theoretic properties of the symmetric integer matrix $A(r)$ to representations of the group $\pi(r)$ into the Lie group  $SU(2)$, should be considered as a Langlands-ian theorem already. Despite its primitivity, nevertheless, it should be clear that this theorem has the potential to disrupt any candidate for a {\it rigorous} $4D$ $\mathcal N$=4 SUSY YM theory, (see \cite{W}, p.12), just as the also very primitive, Rohlin's theorem disrupts low-dimensional handle-body theory, compare to \cite{S}.

 The sheer existence of the above AP Langlands Theory should be an obstruction to rigorizing {\it all} of  \cite{KW} purely analytically, in a similar fashion to the AP-obstructions to the problems discussed in \cite{W}, see also remark 1 in section 6 below.

In is our contention that the discrete purely group-theoretic Artin Equation and the discrete purely group-theoretic, cone-like, {\it $\infty$-generated at each stage}, topology-changing, but homology-preserving, action of the Torelli, with its chaotic finite covering theory (see remark 1 in section 6 below), i.e. our new rigorous, cone-like, tempered dynamic SUSY, \cite{W}, p.6, is not very compatible with classical analytic and Hilbert space using methods. For example, it is hard to imagine a classical Morse Theory, or a classical deformation theory, which group-theoretically changes $4D$-smooth structures , but preserves the underlying topological structure, see \cite{C}, \cite{W},p.5, p.12; see also the knot theory example in \cite{W3}, p.2, \cite{W1}, p.8.

\section{ Essential Differences Between the two Theories}

Although the two theories have some highly non-trivial similarities, e.g., both being related to Donaldson/Seiberg-Witten theory, \cite{Wi}, \cite{W1},p. 9, \cite{CW}, some important general differences are:

i) AP Theory has a rigorous Planckian starting point, see \cite{W}, p.2: the smooth diffeomorphism $h(r):\Omega_n\to \Omega_n$ is determined, up to boundary-fixing smooth isotopies, by the {\it discrete} purely group-theoretic Artin presentation $r$.

 AP Langlands is a mathematically rigorous theory, not a heuristic path integral-using, nor classical SUSY using  theory, e.g., it has a rigorous topological analogue of Montonen-Olive {\it conjectured} $S$-duality.

It is unifyingly conceptually simple and 'complete' mathematically, e.g., {\it the $\sigma$-model picture never breaks down in AP Langlands theory}. There is no need to differentiate 'unramified' from 'ramified' in AP Langlands theory, compare to \cite{GW1}, \cite{F3}.

AP Theory starts from the discrete, i.e., zero-dimensional, Artin presentation $r$, compare to \cite{Wi3}. The highest real dimension in AP Theory is $4D$, and even smooth $4D$ structures are all in function of $2D$ membranes $h(r):\Omega_n\to \Omega_n$, determined by $r$, up to boundary fixing isotopies of the compact $2$-disk, $\Omega_n$, with $n$ holes, see \cite{W}, p.2.

All manifolds in AP Theory are compact and smooth and of real dimension less or equal to $4$.There are, a priori, no analytic moduli spaces in AP Theory.

The only $2D$ manifolds in AP Theory are the compact $2$-disks with $n$ holes, the $\Omega_n$, the only $4D$ manifolds are the compact, simply-connected, smooth  $W^4(r)$, whose always connected, orientable, closed $3D$ boundary, $M^3(r)$, can represent, up to homeomorphism, {\it any} connected, orientable, closed 3-manifold, see, e.g., \cite{W1}, p.8.
 
The class of smooth 4-manifolds considered by AP Theory, the $W^4(r)$, is much more general than that of \cite{KW}; even the {\it closed} ones, i.e.,  when $\pi(r)$ is the trivial group,  contains all complex, simply-connected, elliptic $E(n)$ surfaces, in particular, the Kummer surface; see \cite{CW}.

{\it Thus the $3D$ Poincar\'e Conjecture implies that a very non-trivial class of closed, smooth, simply-connected $4$-manifolds is defined simply by representing them with the Artin presentations of the trivial group.} 

We conjecture that, after the exact relationship of Kapustin-Witten theory and AP theory has been established, these closed, smooth, simply-connected $4D$ manifolds, will eventually  serve, in AP Theory, as the analogues of their higher-dimensional Calabi-Yau manifolds of \cite{KW}.

All analytic, nomographical, classic SUSY difficulties are avoided, since the fundamental Artin Equation is {\it not} a analysis-of-any-kind using, nor SUSY-using equation, as explained on p.2 of \cite{W}.

For example, difficult analytic moduli problems, especially when these are non-compact, do not arise in AP theory, nor is there a need, as already mentioned in the introduction, to {\it 'place some supersymmetric boundary condition..}', as on p.112 of \cite{KW}, \cite{GaWi}.

Instead of a specific classical $S$-duality being the analogue of Langlands duality, as in the \cite{KW} theory, this duality is now represented in AP Theory, by the cone-like, $\infty$-generated at each stage, graded group of homology-preserving, but smooth topology-changing Torelli transitions/interactions, see \cite{W}, p.5, which, furthermore, acts naturally in unison on the $r$, the $h(r):\Omega_n\to \Omega_n$ and the $W^4(r)$. Compare to \cite{Wi}, p.33. They incarnate  at the same time 'mirror' symmetry for our canonical $\sigma$-model, the whole of AP Theory, compare to p. 60 of \cite{Wi}.

This new model-free symmetry, this powerful Torelli action, is also very 'chaotic', due to finite covering theory, see remark 1 in section 6 below.

All this is accomplished so meta-mathematically 'tightly' that a purely discrete group-theoretic analogue of Donaldson/Seiberg-Witten theory is still at hand, \cite{W4}, p.240, \cite{R}, p.621, \cite{CW}, p.2, \cite{W1}, p.9.

ii) Intrinsicity and model-independence.

AP Langlands is model-independent, as Cobordism Theory is, and deals with discrete group-theoretic {\it presentation} theory, not with representation theory of some particular reductive  Lie group $G$ and its Langlands dual. The only Lie group which appears naturally, a priori, in AP Theory, is $SU(2)$, via Donaldson theory, \cite{R}, p.621, \cite{W4}, p.240, \cite{W1}, p.9.

iii) The most extreme mathematically reduction in AP Theory,of the $4D$ $W^4(r)$, first to the $0D$, i.e., discrete, Artin presentation,$r$, then to the to $2D$ $h(r):\Omega_n\to \Omega_n$, i.e., compactification.

 The heuristic {\it "$\Sigma$ large, $C$ small"} argument of \cite{KW}, appears rigorously {\it ab initio} as a 'fait accomplis', so to speak, in AP Langlands. There is no need to 'twist' $4D$ $\mathcal N$=4 Super YM Theory, in a elaborate heuristic manner, to a $2D$ $\sigma$-model, which is a starting point of \cite{KW}.

The $2D$ flat membranic hadrons $h(r):\Omega_n\to \Omega_n$, acted upon by the Torelli action, satisfy Witten's \cite{Wi}, p.56, {\it 'least energy, most SUSY'}, of flat connections on a Riemann surface, in a sharper smooth topological manner.

iv) The intimate relation of AP Langlands with $3$-manifold theory, (instead of complex Riemann surfaces), including the complete Knot and Linking theory in any closed, orientable $3$-manifold, see remark 5 in section 5 ahead.

Any closed, orientable $3$-manifold, and simultaneously {\it any} link therein, not just its fundamental group, can be represented by an Artin presentation.

There is an immediate clear relation, in AP Theory to Knot and Linking Theory in any closed $3$-manifolds, see\cite{W4}, p.227, \cite{W3}, \cite{W1}, p.8; compare to \cite{Wi}, p.5, \cite{Wi8}.

 Since the Thurston Geometrization Program is obviously a $3D$ Automorphic program, AP Langlands has again the  right to call itself 'geometric'. This is augmented by the Nielsen-Thurston theory of the diffeomorphisms $h(r):\Omega_n\to \Omega_n$; see remark 4 in section 5 ahead and eventually \cite{W2}.

v) A closer immediate relation with L-functions. 

 The first zeta-functions appearing in our Langlands Theory are the zeta-functions of Artin-Mazur, for the diffeomorphisms $h(r):\Omega_n\to \Omega_n$, \cite{AM}, \cite{Ru}, which, \'a la Langlands, should 'match', at least be related to, the Hecke, Maass $\theta$- and L- functions, (see \cite{M}), pertaining to the symmetric integer matrices $A(r)$ of p.7,\cite{W1}; see also remark 1 in  section 5, ahead.

{\it This is a very important Langlandsian property which does not appear very clearly in \cite{KW}}.

 If it did, it would probably be related to heuristic path integral correlations.

vi) More clarity and unification in AP Operator Theory.

 In AP Langlands, since it is parameter-free, all operators are 'rigid' in the sense of \cite{GW}, and they are all part of the same graded group and furthermore act in unison: as string operators on the Artin presentations $r$, as surface operators on the $h(r):\Omega_n\to \Omega_n$ and as Hecke operators on the $4D$ $W^4(r)$. Compare to \cite{GW}, p.71, \cite{KS}, \cite{KSV}. 

The action of these Torelli operators can be very subtle with respect to Knot Theory, see the example on p.2 of \cite{W3}, p.8 of \cite{W1}.

 We conjecture that, via finite covering theory, 'operator theory' becomes too chaotic to be studied just with analytic and other classical methods, see remark 1 in section 6, below.

vii) An immediate relation to analogues of other heuristic theories of Modern Physics, e.g., QCD and confinement, holography, string theory, dynamic dark energy, YM Millenium Problem and Mass gap, etc., as explained in \cite{W}, \cite{W1}.

  AP Theory is, in particular, a rigorous topological analogue of the Strong Force, see \cite{W},p.5, compare to \cite{H}, \cite{HSZ},\cite{Wilc}, augmented {\it dynamically} by by the cone-like, $\infty$-generated at each stage, graded topology-changing, but homology-preserving Torelli operators.
 This is the qualitative 'membranic' analogue of Montonen-Olive's EM theory in $4D$ $\mathcal N$=4 Super YM theory.

 We conjecture that AP Langlands theory, at the purely $2D$ level, contains deep in itself, an ultimate Confinement Theory, an AP-analogue to Nielsen-Thurston theory where the 'quarks' and 'partons' as fixed and periodic points of the $h(r):\Omega_n\to \Omega_n$, are 'confined' into the Artin presentations $r$ as explained in \cite{W}, p.4, section 5.

viii) AP Langlands theory is {\it complete} as a $\sigma$-model and as discrete group-theoretic, i.e., model-free gauge-theoretic as possible. AP Theory contains as many 'bonus' symmetries as possible, compare to \cite{I}.

  This is very relevant to the DARPA's Math challenge no. 17, "Geometric Langlands and Quantum Physics": "How does the Langlands Program, which originated in number theory and representation theory, explain the fundamental symmetries of Physics?".

\section{More Comments on the Above}

 Since AP Theory is a somewhat unconventional, non-analytic, non-categorical, very elementary and conceptually simple, but subtle theory, at the risk of excessive repetition, we make some  more pertinent general comments which are relevant to and might help clarify some of the above.

The Artin Equation, \cite{W1}, p.7,  \cite{W}, p.2, is a discrete purely group-theoretic equation for the Artin presentations $r$, which encodes, purely  group-theoretically, the $2D$ 'membranic flatness' of the $h(r):\Omega_n\to \Omega_n$. This is the main reason AP-Langlands Theory is so conceptually simple, since no analytic, nor nomographical difficulties can occur. Nor does it use any 'sewing', 'fusion' rules nor 'bordism glueing' of any kind, despite its radical topology-changing Torelli dynamics.

Artin presentations are very important purely group-theoretically, because they characterize the fundamental groups of closed, orientable 3-manifolds, \cite{W1}, p.6; i.e., an arbitrary abstract group is isomorphic to the fundamental group of a closed, orientable $3$-manifold, if and only if it has an Artin presentation.  These groups now form the intrinsic AP analogue of Loop Groups.

 AP Langlands theory, this ultimate meta-mathematical smooth topological reduction of Geometric Langlands Theory, containing non-trivial analogues to all of the 'main ideas' of the introduction to \cite{KW}, is natural, due to the radical canonical compactification III. of section 2 above, giving a canonical $2D$ $\sigma$-model, targeting the smooth $4D$ manifolds $W^4(r)$, which together with the Torelli action, incarnates the ultimate quantum topological Langlands theory of this paper. Perhaps it should be called Membranic Langlands theory, as opposed  to  Topological/Geometric Langlands theory.   

This model-independent, universal $\sigma$-model, i.e. all of AP Theory, is where AP-Langlands Theory lives and, a priori, there is no need for Category Theory, classical heuristic SUSY methods, etc., in order to describe genuine qualitative topological analogies to those of \cite{KW}. See also remark 4 in section 6 below.

This chamaleonic 'transfer' behaviour of AP-Theory is not all that surprising due to its basic group-theoretic Erlanger Program -like nature. It consists of pure discrete, topology-less, group theory, augmented topologically only by adjoining the smooth topological diffeomorphisms, the membranic hadrons  $h(r):\Omega_n\to \Omega_n$, which then define the smooth $4D$-branes, the $W^4(r)$, where the Torelli transitions operate in a $Z$-homology preserving, but smooth topology changing manner. Compare to \cite{K2}.
It consists of the maximum of symmetry adjoined by the minimum of smooth topology, i.e., $2D$ flat membranicity.

We can say Intrinsic Quantum Topological Langlands Theory lives in a huge, super $2D$ $\sigma$-model targeting compact, smooth, simply-connected $4D$ manifolds, with a connected boundary, namely AP Theory itself. The $\sigma$-model picture is everything; it completely encompasses our enveloping  completion, compactification of $4D$ $\mathcal N$=4 SUSY YM Theory.

With this radical compactification, AP Theory by its sheer existence gives a genuine intrinsic (i.e., model-free) quantum topological Langlands Theory, avoiding classical SUSY (as explained on p.6 of \cite{W}), such difficult to handle analytic concepts as non-compact moduli and such esoteric category-theoretic concepts as sheaves, eigensheaves,  derived categories, $\mathcal D$-modules, stacks, etc.,and, furthermore, embedding all 'operators' into one single graded discrete group, where now they should be studied, not just algebraically, but {\it geometrically}, in a genuine Erlanger Program style, {\it within} the graded AP Torelli group, compare to \cite{KS}.  See also remark 1 in section 6 below.
 
 Among all the Langlands theories, AP Langlands theory is the compact, universal topological envelope theory of \cite{KW}'s Geometric Langlands Theory from the point of group theory, i.e.,intrinsic, non-perturbative, model-free discrete gauge theory, compare to \cite{Wi4}, p.1. In fact,in AP Theory, the Geometric Langlands Program actually  becomes a Super Erlanger Program for Physics.

Perhaps the ultimate Geometric Langlands correspondence should consist of relating these two theories more intimately, starting by relating \cite{KW}'s classic, static sense of SUSY to the new, more flexible, more {\it 'dynamic'}, rigorous, cone-like, $\infty$-generated at each stage, graded, membranically {\it tempered} 'SUSY' described above and in \cite{W}, p.6. All this in a time independent manner, etc., see \cite{KW}, p.123.

 AP Langlands is more than a non-perturbative Gauge Theory, it has  a {\it new} type of cone-like,  $\infty$-generated at each stage, graded, {\it model-independent symmetry}, but still containing a non-trivial purely discrete group-theoretic Donaldson/Seiberg-Witten theory, \cite{W1}, p.9, \cite{CW}.

 Thus AP Langlands Theory is indeed a universal, model-independent Erlanger Program for Modern Physics.

\section{New Beacons in AP Theory for Langlands Theory?}

1.  AP Langlands suggests very strongly that the Artin-Mazur zeta-functions, \cite{AM}, \cite{Ru}, of the diffeomorphisms $h(r):\Omega_n\to \Omega_n$ should be related to the L-functions of Siegel Modular Theory, as in, say, Maass, \cite{M}, pertaining to the integer symmetric matrices, $A(r)$, which determine the $Z$-homology and quadratic form of $W^4(r)$, as well as the $Z$-homology of its connected boundary $M^3(r)$. 

Physically, this would be the ultimate Langlands Functoriality relation: L-functions related to the {\it vacuum} fluctuations, the hadrons $h(r):\Omega_n\to \Omega_n$, our most {\it local} topological concept, would be related to L-functions pertaining to the spacetime {\it universes} $W^4(r)$, our most {\it global} topological concept. Compare to Langlands': {\it "Functoriality is the core notion of what is frequently referred to as the Langlands program"}.\cite{L1},p.3.

Should this, in unison with the 'A-branes' of \cite{KW}, help in clearing up some of the mysteries regarding the role of the {\it symplectic} group in Siegel's theory of quadratic forms, as described by Weil in \cite{We}?

An important starting point here is to investigate what occurs under the action of the Torelli, since they always preserve the symmetric integer $n\times n$ matrices $A(r)$.

2. An intrigueing feature from AP theory is the 'double presentation problem': the Artin presentation $r$, since it also defines a connected, closed, orientable $3D$-manifold, $M^3(r)$, namely the boundary of $W^4(r)$, it also  determines the mapping class group of the former, call it $\Gamma(r)$. It is known to be a finitely presented group, very important in classical $(3+1)$-TQFTs, \'a la Atiyah, which according to the Thurston Geometrization Program should be related to $\pi(r)$ in a Nielsen-like simple manner.

Thus, two problems arise:
 i) from $r$ determine $\Gamma(r)$ purely group-theoretically, e.g., find a presentation of it, given the Artin presentation $r$. 
ii) Is this problem an AP analogue, a surviving vestige under model-independence, of some Geometric Langlands problem? Perhaps, given $r$,  $\Gamma(r)$ is some {\it 'endoscopic'} group of $\pi(r)$?

3. Let $U^4(r)$ be the compact smooth $4$-manifold constructed on p.622 of \cite{R}, where the diffeomorphism there is the one obtained from $h(r)$ by extension by the identity, to the boundary of $\Omega_n\times I$, with I = [0,1] and $\Omega_n$ identified with $\Omega_n\times 0$; then, as a result of the planarity of the $h(r):\Omega_n\to \Omega_n$, it is easy to see that its fundamental group is isomorphic to the {\it direct}, (not free) product $\pi(r)\times Z$.

{\it Since all of the AP Langlands theory above followed from our fundamental $4D$-brane construction of the $W^4(r)$,from the $r$,  via the $h(r):\Omega_n\to \Omega_n$, \cite{W1}, p.4, what Langlandsian analogy should follow using $U^4(r)$ instead of $W^4(r)$?} 

Is this related to any Riemann-Hilbert correspondence in classical Geometric Langlands Theory, \'a la Laumon-Drinfel'd? Compare to \cite{K1}and \cite{F6}.

To any Flux Theory in \cite{KW}? Compare to \cite{K}.

 What does it mean Langlands-wise that the natural infinite cyclic cover of $U^4(r)$, corresponding to the $Z$ factor, has a fundamental group isomorphic to $\pi(r)$, and boundary consisting of an infinite disjoint union of copies of $M^3(r)$, each again having $\pi(r)$ as a fundamental group, and one component homeomorphic to $\Sigma_n\times R$, where $\Sigma_n$  is the boundary of $\Omega_n\times I$?

4. In \cite{G}, Gonz\'alez-Acu\~na, using AP Theory, gave a purely discrete group-theoretic equivalent of the Poincar\'e Conjecture involving, at each stage $n$, the existence of a certain {\it double coset} in the automorphism group of a certain discrete finitely presented $1$-relator group. 

Does this have an analogue to any {\it automorphic} Hecke sheaf arguments in classical Geometric Langlands theory?

5.  Geometric Langlands theory and Knot and Linking theory?

AP Langlands Theory is intimately related to the complete Knot and Linking theory in any closed, orientable $3$-manifold, $M^3(r)$, \cite{W4}, p.226, \cite{W3}, \cite{W1}, p.8, where the main knot and linking invariants  are easily obtainable from the Artin presentation $r$, via the algebraic computer system MAGMA. Compare to \cite{Wi}, p.5, and \cite{Wi8}.

6. If we consider the Artin presentation $r$ as a 'quantum' string, see \cite{W}, p.7, is it physically relevant, due to its strong Torelli dynamics, to Langlands' Planar Percolation theory, e.g., does it give, e.g.,  a 'limit form' of planar percolation?  Does it establish 'universality' in percolation?, see \cite{L2}, p.9.

\section{Remarks, Problems, Questions, Conjectures} 

1. Due to Perelman's proof of Thurston's Geometrization Conjecture, we now know that the groups $\pi(r)$ are all {\it residually finite}. This implies that AP Theory has a very rich finite regular covering theory for the $M^3(r)$, each new fundamental group again having an Artin presentation.

 When combined with the Torelli action, as well as the iteration of the $h(r):\Omega_n\to \Omega_n$, (the AP Frobenius map), this leads to a very  {\it chaotic} situation in AP Langlands' operator theory: e.g., despite that the universes $W^4(r)$ are simply-connected, the whole set of them still get acted on by the many normal subgroups of finite index of the fundamental group, $\pi(r)$, of their connected boundary, $M^3(r)$.

 What does this really  mean for \cite{KW}, classic Geometric Langlands Theory, Modern Physics and Cosmology? Compare to \cite{BKLS}, p.14.

2. What is the analogue in AP Theory of \cite{Ge}, p.178: {\it "That the possible number fields of degree $n$ are restricted in nature by the irreducible $\infty$-dimensional representations of $GL(n$), was the visionary conjecture of R. P. Langlands."} ?

What will be restricted in this manner by the cone-like, {\it $\infty$-generated at each stage}, graded group of topology-changing, but homology preserving Torelli transitions/interactions above? 

Just the topology of the $M^3(r)$? Is this the ultimate cause of the truth of the Geometrization Conjecture?  A priori, in AP Theory, there was no reason at all to conjecture that any closed, orientable 3-manifold $M^3(r)$, where $\pi(r)$ is trivial, should be homeomorphic to the 3-sphere, $S^3$, compare to \cite{W4}.

If we substitute this $\infty$ for the 'unitarian' one used by Weil, \cite{We}, mentioned in the previous section, do we shed more light on the role played by the {\it symplectic group}, via the symmetric integer matrices $A(r)$, on the celebrated work of C. L. Siegel on quadratic forms? \cite{We},\cite{M}.

3. Is there a general meta-mathematical procedure to go, in general, (call it 'AP-ization') from certain heuristic classical SUSY using theories (e.g., $4D$ $\mathcal N$=4 Super YM theory, Kapustin-Witten theory,..,) to rigorous, {\it model-free}, dynamic AP-tempered SUSY using theories (e.g., AP $(3+1)$-QFT, AP Langlands theory,..,), while still preserving crucial important analogies? This would lead to new physically generated mathematically rigorous theories.

4. Is category theory really ultimately necessary and fundamental in the classic Geometric Langlands Program?
   Does the sheer meta-mathematical existence of the above AP Langlands Theory, the above super Erlanger Program, imply a negative answer? 

In other words, does the embedding of  AP Geometric Langlands into the canonical universal $\sigma$-model, that is AP Theory, {\it lock} non-trivial Category Theory {\it out}, so to speak, due to AP theory's topological Planckian starting point and its strong, chaotic, $\infty$-generated at each stage graded, topology-changing, but $Z$-homology preserving Torelli action, combined with finite covering theory (as in remark 1 above), as well as the presence still of a purely group-theoretic Donaldson/Seiberg-Witten Theory?  The same should occur to classical analytic differential-geometric and classical SUSY and Hilbert space using methods, as stressed in \cite{W}?

We conjecture that the answer to the last two questions is in the affirmative.

\end{document}